\documentclass[12pt,a4paper]{article} 
\usepackage{fullpage,amsmath,amssymb,amsthm}

\DeclareMathOperator{\Pol}{Pol}

\newtheorem{lemma}{Lemma}[section]
\newtheorem{theorem}[lemma]{Theorem}

\theoremstyle{definition}
\newtheorem{defn}[lemma]{Definition}
\newtheorem{rem}[lemma]{Remark}
\newtheorem{sub}[lemma]{}

\begin{document}

\title{\textbf{Polynomial Functions on a Central Relation}}

\author{Dietmar Schweigert}

\date{\today}
\maketitle

\thispagestyle{empty}

\setcounter{page}{0}
\vspace{2cm}

\abstract{We show that the algebra $\mathbf{R} = \bigl(R; \wedge,
  \veebar, 0, \bar{f}_i(x) (i \in I)\bigr)$ is central polynomially complete
  and finite.}

\bigskip

\noindent\begin{small}Keywords:  Central polynomially complete, prepolynomially
  complete: AMS Mathematical Subject Classification 08A40, 06F99.
\end{small}

\newpage

\section{Introduction}

A clone $C$ on a set $X$ is a set of finitary operation $f : X^n \to X$, which
contains all the projections and is closed under composition.  We can also
consider a clone as a set of term function of some universal algebra.  We have
to choose a universal algebra $\mathbf{A}$, which generates the clone $C$.

We consider clones of binary central relations and choose an algebra
$\mathbf{R} = \bigl(R; \wedge, \veebar, 0$, $\tilde{f}_i (i \in
I)\bigr)$ which preserves the central relation.   $\mathbf{R}$
consists of lattice operations without $0$ and contains unary operations.   We
shall show that the algebra is central polynomially complete and finite.   As
the clone of the central relation is maximal we discuss two problems on
prepolynomial algebras.   Our approach is similar to the theory of order
polynomially complete \cite{Schw2}, \cite{PH}, \cite{GS}, \cite{Go1}.

\section{Prepolynomially Complete Algebras}

Let $A$ be a set and $f : A^n \to A$ a function in $n$ variables (that is, an
operation).   A clone $H$ is a set of operations on $A$, which is closed by
composition of functions and by its manipulation on variables (for the
definition see \cite{Ro1}, \cite{Ro2}).

Let $\rho$ be an $h$--ary relation on $A$.   Then $\Pol_\rho$ is the clone of
function which preserves the relation $\rho$.

Let $\mathbf{A} = (A; \Omega)$ be an algebra of type $\tau$.   The
clone $P(A)$ of polynomial functions of $\mathbf{A}$ consists of all
$n$--place polynomial functions $n \in \mathbb{N}$.   We define recursively

\begin{enumerate}
\item The projection $e_i^n : A^n \to A$ with $e^n_i (x_1, \dots, x_n) = x_i$
  and the constant function $c^n_a : A^n \to A$ with $c_a^n(x_1, \dots, x_n) =
  a$, $a \in A$ are $n$--place polynomial functions.

\item If $f \in \Omega$ is an $m$--place operation of the algebra
  $\mathbf{A} = (A; \Omega)$ and $p_1, \dots, p_m$ are $n$--place
  polynomial functions, then $f\bigl(p_1(x_1, \dots, x_n), \dots, p_m(x_1,
  \dots, x_n)\bigr)$ is a polynomial function.
\end{enumerate}

One may apply (2) only a finite number of times.

\begin{defn}\label{defn21}
  Let $\mathbf{A} = (A; \Omega)$ be an algebra.  A relation $\rho \in
  A^n$ is $\rho$--compatible if it holds for every operation $f \in \Omega$.
  If $(a_{11}, \dots, a_{h1}) \in \rho, \dots, (a_{1n}, \dots, a_{hn}) \in
  \rho$, then $\bigl(f(a_{11}, \dots, a_{1n}), \dots, f(a_{h1}, \dots,
  a_{hn})\bigr) \in \rho$.
\end{defn}

\begin{rem}\label{rem22}
  Every polynomial function of the algebra $\mathbf{A}$ is compatible.
\end{rem}

\begin{defn}\label{defn23}
  An $n$--ary relation $\rho$ is central if there is a non--empty proper
  subset $Z$ of $A$ such that

  \begin{enumerate}
  \item $(a_1, \dots, a_n) \in \rho$ if at least one $a_i \in Z$;

  \item $\rho$ is invariant under permutations of coordinates;

  \item $(a_1, \dots, a_n) \in \rho$ if $a_i = a_j$ for some distinct $i,j$.
  \end{enumerate}

The subset $Z$ of $A$ is called a center.   We shall consider only binary
central relations on a set $A$ with the center $\{0\}$ for some fixed element
$0 \in A$.
\end{defn}

\begin{defn}\label{defn24}

  \begin{enumerate}
  \item A central function $f : A^n \to A$ is a function which is compatible
  on a central relation $\rho$.

\item A central polynomial $f : A^n \to A$ is a central function which is a
  polynomial function.
  \end{enumerate}
\end{defn}

\begin{defn}\label{defn25}
  The algebra $\mathbf{A} = (A; \Omega)$ is polynomially complete
  (functionally complete) if every function $f : A^n \to A$ is a polynomial
  function for $n \in \mathbb{N}$.

Examples can be found in the book of Lausch and N\"obauer \cite{LN}.
\end{defn}

\begin{defn}\label{defn26}
  The algebra $\mathbf{A} = (A, \Omega)$ is prepolynomially complete
  if every function $f : A^n \to A$ is a polynomial function of the algebra
  $\mathbf{A} = (A; \Omega \cup \{g\})$ for every $g \not\in P(A)$.
\end{defn}

\begin{defn}\label{defn27}
  \begin{enumerate}
  \item The algebra \mbox{$\mathbf{A} = (A, \Omega)$} is
    $n$--central--polynomially complete, if every $n$--ary function \mbox{$f :
      A^n \to A$,} $n \in \mathbb{N}$ is a polynomial function on the central
    relation $\rho$ of $A$.

\item The algebra \mbox{$\mathbf{A} = (A, \Omega)$} is central--polynomially
  complete, if every function\\ \mbox{$f : A^n \to A$} for every place $n$ is a
  polynomial function on the central relation of $A$.
  \end{enumerate}
\end{defn}

\section{$\rho$--Polynomially Complete Algebras}

\begin{defn}\label{defn31}
  Let $\rho$ be an $n$--place relation on the finite or infinite set $A$.
  The algebra $\mathbf{A} = (A, \Omega)$ is $n$--$\rho$--polynomially
  complete if every $\rho$--compatible function in $n$ variables on $A$ is a
  polynomial function in $n$--variable.

The algebra $\mathbf{A}$ is $\rho$--polynomially complete if every
$\rho$ compatible function on $A$ is a polynomial function for every variable,
that is, $P(A) = \Pol_\rho$.
\end{defn}

\begin{theorem}\label{thm32}
  If an algebra $\mathbf{A} = (A, \Omega)$ is
  $n$--$\rho$--polynomially complete for every $n$, then the algebra
  $\mathbf{A}$ is also 1--$\rho$--polynomially complete.
\end{theorem}

\begin{proof}
  Let $f$ be a $\rho$--compatible function in one variable $x_n$.   We define
  a function $g$ for all $(x_1, \dots, x_n) \in A^n$ by $g(x_1, \dots, x_n) =
  f(x_n)$.   Because $f$ is $\rho$--compatible, then the function $g$ is also
  $\rho$--compatible.   By assumption, $g$ is a polynomial function.

Let $a \in A$.   We put $x_1 = a, \dots, x_{n-1} = a$, then we obtain a
polynomial function $g(a, \dots, a, x_n)$ in one variable.
\end{proof}

\section{A Central--Polynomially Complete Algebra}

\begin{sub} We consider a finite set with an element $0 \in R$ and we
consider a binary relation on $R$
\[
\rho : \bigl\{(r,0) \mid r \in R\bigr\}\; \cup\; \bigl\{(0,s) \mid s \in
R\bigr\}\; \cup\; \bigl\{(t,t) \mid t \in R\bigr\}\,.
\]
$\rho$ is reflexive, symmetric and is a central relation.   $\rho$ is not
transitive.   In the case that $\rho$ would be transitive, then it follows
from $r_1 \rho\, 0$ and $0\, \rho\, r_2$ to $r_1 \rho\, r_2$ and, therefore,
$\rho$ would be the all relation.
\end{sub}

\begin{sub} We consider the algebra
\[
\mathbf{R} = (R; \wedge, \veebar, 0)
\]
where $(R; \wedge, \veebar, 0)$ is a semi--lattice with a zero element $0$
with $0 \wedge x = x \wedge 0 = 0$.

The operation $\veebar$ is defined by
\[
x \veebar y = 
\begin{cases}
  x \vee y & \text{ if $x \not= 0$ and $y \not= 0$}\\
0 & \text{ else}.
\end{cases}
\]
\end{sub}

\begin{rem}\label{rem43}
  The operations $\wedge$ and $\veebar$ are compatible.

From the definition of $\rho$ it follows that $a = b$ for every $a \not= 0$
and $b \not= 0$.   We have $a \wedge 0 = 0 \wedge b = 0$.

If $(a_1, b_1) \in \rho$ and $(a_2, b_2) \in \rho$, then we have
$\bigl(f(a_1, a_2), f(b_1, b_2)\bigr) = (a_1 \wedge a_2, b_1 \wedge b_2) \in
\rho$.

If $(a_1, b_1) \in \rho$ and $(a_2, b_2) \in \rho$, then

\begin{enumerate}
\item $(a_1 \vee a_2, b_1 \vee b_2) \in \rho$ in the case $a_1 \not= 0$,
  $a_2 \not= 0$, $b_1 \not= 0$, $b_2 \not= 0$, which means: $(a_1 \vee a_2,
  a_1 \vee a_2) \in \rho$.

\item $(a_1 \veebar a_2, b_1 \veebar b_2) = (0,0) \in \rho$.
\end{enumerate}
\end{rem}

\begin{sub}\label{sub44}
  We consider all one--place functions $f_i : R \smallsetminus \{0\} \to
  R$ and define all the one--place functions $\bar{f} : R \to R$
\[
\bar{f}_i(x) = 
\begin{cases}
  f_i(x) & x \not= 0\\
0\,. &
\end{cases}
\]
Note that $\bar{f}_i$, $i \in I$, are compatible.

We have the algebra
\[
\mathbf{R} = \bigl(R; \wedge, \veebar, 0, \bar{f}_i (i \in I)\bigr)\,.
\]
\end{sub}

\begin{sub}\label{sub45}
  The algebra $\mathbf{R}$ is 1--$\rho$--polynomially complete because
  it contains all compatible unary functions.
\end{sub}

\begin{sub}\label{sub46}
  The algebra $\mathbf{R}$ is $n$--central polynomially complete.
\end{sub}

\begin{proof}
  We consider compatible functions in $n$ variables $n > 2$
\[
f(x_1, \dots, x_n) = \veebar f_{(a_1, \dots, a_n)} (x_1, \dots, x_n)\,,
\]
where 
\[
f_{(a_1, \dots, a_n)} (x_1, \dots, x_n) = 
\begin{cases}
f(a_1, \dots, a_n) & x_1 = a_1, \dots, x_n = a_n\\
0 & \text{else}
\end{cases}
\]
and where $f_{(a_1, \dots, a_n)} (x_1, \dots, x_n) = 0$ if $a_1 = 0$ or
$\ldots a_n = 0$.

Clearly the function $f_{(a_1, \dots, a_n)} (x_1, \dots, x_n)$ is compatible.

We like to decompose more
\[
f_{(a_1, \dots, a_n)} (x_1, \dots, x_n) = g_{(a_1, \dots, a_{n-1})} (x_1,
\dots, x_{n-1}) \veebar h_{(a_n)} (x_n)
\]
where
\[
g_{(a_1, \dots, a_{n-1})} (x_1, \dots, x_{n-1}) = 
\begin{cases}
  f(a_1, \dots, a_n) & x_1 = a_1, \dots, x_{n-1} = a_{n-1}\\
0 & \text{else}
\end{cases}
\]
and where
\[
h_{(a_n)} (x_n) = 
\begin{cases}
  f(a_1, \dots, a_n) & x_n = a_n\\
0 & \text{else}\,.
\end{cases}
\]
By our assumption $\mathbf{R}$ is 1--$\rho$--polynomially complete
and, therefore, $n$--$\rho$--polynomially complete.
\end{proof}

\begin{rem}\label{rem47}
  For our example of the central complete algebra we have used the same
  construction of the partial clone as in \cite{HRS}.
\end{rem}

\section{The Cardinality of the Central Polynomial Functions}

\begin{theorem}\label{them51}
  Let $(A; \Omega)$ be an infinite algebra which preserves a central relation
  $\rho$ with the center $\{Z\}$.  Let the type of the algebra be finite
  (which means $|\Omega|$ is finite).  Then the algebra $(A; \Omega)$ is not
  1--central polynomially complete.
\end{theorem}

\begin{proof}
  We count the set of the polynomial functions on the algebra
  $\mathbf{A}$ and assume that $A \smallsetminus Z$ is infinite..
  
  Let $W = W\bigl(A \cup \{x\}\bigr)$ be the word algebra over $\Omega$.  Let
  $C_n$ be the class consisting of all words of length $n$.  $C_n$ is a subset
  of $\bigl(\Omega \cup A \cup \{x\}\bigr)^n$, hence $|C_n| \le
  \big|\bigl(\Omega \cup A \cup \{x\}\bigr)\big|^n = |\Omega \cup A \cup
  \{x\}| = |\Omega| \cup |A| \cup |\{x\}| = |\infty| \cup  |A| \cup 1 = |A|$.

Since $W = \cup (C_n \mid n \ge 1)$ we have $|W| = \Sigma \bigl(|C_n|\; n \ge
1\bigr) \le \Sigma \bigl(|A|\, h \ge 1\bigr) = |A|$.   As every element of
$P_1(A)$ can be presented as a word of $W$, there is injection from $P_1(A)$
to $W$, hence $\big|P_1(A)\big| \le |W| \le |A|$.

Now we have to count the set of the central functions in one variable.  We
estimate $A \smallsetminus \{Z\}$ by $|A\smallsetminus Z|^{|A\smallsetminus
  Z|}$ because it contains all functions from $A \smallsetminus \{Z\}$ to $A
\smallsetminus \{Z\}$. In the case that $A \smallsetminus Z$ is finite and the
center $Z$ is infinite, we have the similar argument.
\end{proof}

\section{Patterns of Prepolynomially Complete\\ Algebras}

\begin{sub}\label{sub61}
  An algebra $\mathbf{A}$ is polynomially equivalent to the algebra
  $\mathbf{B}$ if there exists an algebra $\mathbf{A}^\prime$, which is
  isomorphic to $\mathbf{A}$ such that $P(\mathbf{A}^\prime) =
  P(\mathbf{B})$.   We have considered the algebra $\mathbf{R} = (R;
  \wedge, \veebar, 0)$, which seemed to be special.   Of course, we could also
  choose another algebra, which is polynomial equivalent to $\mathbf{R}$.
  The choice of the algebra we shall call \emph{pattern}.

\noindent In this way, we should like to present two problems:

\noindent\textbf{Problem 1:} Is every prepolynomially complete algebra finite?

\noindent\textbf{Problem 2:} Give patterns and criteria of prepolynomially
complete algebras!
\end{sub}

\begin{sub}\label{sub62}
  A clone $C$ of functions on set $A$ is called maximal clone if, for any $f
  \in O_A \smallsetminus C$ follows $\langle (C \cup \{f\})\rangle = O_A$.
  That is, the set $C \cup \{f\}$ generates the clone of all functions on $A$.
  It is clear that maximal clone presents a prepolynomially complete algebra
  on finite set $A$.
\end{sub}

\begin{sub}\label{sub63}
  On the other hand, we may consider the clone generated by the projections.
  A pattern of an algebra is the algebra $\mathbf{A} = (A; \circ)$ with the
  axioms $x \circ x = x$, $(x \circ y ) \circ z = x \circ (y \circ z)$, $x
  \circ y \circ z = x \circ z$.

\noindent Let $A_\lambda$ be a set of the cardinality $\lambda$ and denote the
algebra by $\mathbf{A}_\lambda = (A_\lambda; \circ)$.

\noindent Let $R$ be the set of the diagonal relations on $A_\lambda$.   Then
the algebra $\mathbf{A}_\lambda = (A_\lambda. \circ)$ is
$R$--polynomially complete for every cardinality $\lambda$.
\end{sub}

\vfill

\noindent Technische Universit\"at Kaiserslautern\\
FB Mathematik, D--67663 Kaiserslautern\\
Germany\\[1.0ex]
\noindent E--mail address:  Dietmar.Schweigert@web.de

\end{document}